\documentclass[12pt]{article}
\usepackage{amsfonts}
\usepackage{amssymb}
\makeatletter

\@addtoreset{equation}{section}
\def\section{\@startsection {section}{1}{\z@}{-3.5ex plus -1ex minus
 -.2ex}{2.3ex plus .2ex}{\large\bf}}
\def\subsection{\@startsection{subsection}{2}{\z@}{-3.25ex plus%
 -1ex minus -.2ex}{1.5ex plus .2ex}{\sc}}
 
\advance \voffset by -0.9in
\advance \hoffset by -0.8in
\newcommand{\ZZ}{\mathbb{Z}}
\newcommand{\QQ}{\mathbb{Q}}

\newcommand{\CC}{\mathbb{C}}
\newcommand{\tr}{{\rm tr}}
\def\bea{\begin{eqnarray}}
\def\eea{\end{eqnarray}}
\def\halmos{\hbox{\vrule height0.31cm width0.01cm\vbox{\hrule height
 0.01cm width0.3cm \vskip0.29cm \hrule height 0.01cm width0.3cm}\vrule
 height0.31cm width 0.01cm}}
\def\hhalmos{{\unskip\nobreak\hfil\penalty50
        \quad\vadjust{}\nobreak\hfil\halmos
        \parfillskip=0pt\finalhyphendemerits=0\par}}

\def\mLP{\medbreak\noindent}

\def\al{\alpha}
\def\be{\beta}
\def\ga{\gamma}
\def\de{\delta}
\def\ep{\varepsilon}
\def\ka{\kappa}
\def\De{\Delta}
\def\id{{\rm id}}
\let\lan=\langle
\let\ran=\rangle
\let\ten=\otimes
\let\wb=\overline

\newtheorem{theorem}{Theorem}[section]
\newtheorem{lemma}[theorem]{Lemma}
\newtheorem{corollary}[theorem]{Corollary}

\begin{document}
\baselineskip 18pt
\parskip 7pt
\begin{flushright}
ITFA-99-6,  UvA-Math-99-08 \\
math.QA/9904029
\end{flushright}
\vspace{0.7cm}

\begin{center}
\baselineskip 24 pt
{\LARGE Fourier transform and   the Verlinde formula for} \\
 {\LARGE the quantum double of a finite group}

\vspace{1cm}
\baselineskip 18pt
{\Large 
T.H. Koornwinder\footnote{e-mail: {\tt thk@wins.uva.nl}} }\\
KdV Institute for  Mathematics, University of Amsterdam\\
Plantage Muidergracht 24, 1018 TV Amsterdam\\
The Netherlands\\
\vspace{1cm}
{\Large B.J. Schroers \footnote{e-mail: {\tt bernd@maths.ed.ac.uk}} }\\
Department of  Mathematics and Statistics, University of Edinburgh\\
King's Buildings, Mayfield Road, Edinburgh EH9 3JZ\\
United Kingdom\\
\vspace{1cm}
{\Large 
J.K. Slingerland\footnote{e-mail: {\tt slinger@wins.uva.nl}} and 
F.A. Bais\footnote{e-mail: {\tt bais@wins.uva.nl}} }\\
Institute for Theoretical Physics, University of Amsterdam\\
Valckenierstraat 65, 1018 XE Amsterdam\\
The Netherlands

\vspace{0.5cm}

4 October 1999

\vspace{0.3cm}

\end{center}


\date{\today}

\begin{abstract}
\noindent We define a  Fourier transform $S$  for the quantum double 
$D(G)$ of a finite group $G$. Acting on characters of $D(G)$, $S$ and the 
central ribbon element of $D(G)$ generate a unitary matrix
representation of the  group $SL(2,\ZZ)$. The characters form 
a ring over the integers under both the algebra multiplication 
and its dual, with the latter encoding the fusion rules of $D(G)$.
The Fourier transform  relates the two ring structures. We use
this to give a particularly short proof of the Verlinde formula
for the fusion coefficients.
\end{abstract}
\vspace{1cm}

\section{Introduction}

The quantum double $D(G)$ of a 
finite  group $G$ is a quasi-triangular ribbon Hopf
algebra \cite{qugroups} constructed, 
via Drinfeld's double construction \cite{Drinfeld}, 
out of the Hopf algebra $C(G)$ of  ${\CC}$-valued functions on $G$. 
Such quantum doubles arise in physics in 
orbifold conformal field theories \cite{DVVV} and in the classification
of flux-charge composites in the massive phases of (2+1)-dimensional
gauge theories \cite{BDP},\cite{BP}. The mathematical structure of $D(G)$
was clarified in \cite{DPR}. 
There and in \cite{DVVV} it was  also pointed out that 
the set of irreducible representations of $D(G)$ carries a 
representation of the group $SL(2,\ZZ)$ by unitary, symmetric 
matrices. In particular one has a symmetric and unitary matrix
$S$ and a diagonal, unitary matrix $T$ acting on the 
set of irreducible representations which satisfy the modular 
relation $(ST)^3=S^2$ and $S^4=1$. Although perhaps surprising from the 
point of view of Hopf algebras, the appearance of the 
$SL(2,\ZZ)$ action in the representation theory of $D(G)$
is physically motivated by  application
of $D(G)$ in orbifold conformal field theories.
In particular,  it was already pointed out
in \cite{DVVV} that the matrix $S$ plays the role of the Verlinde 
matrix which diagonalises the fusion rules in  orbifold conformal field
theory (for a general  review 
of fusion rules in conformal field theory 
we refer the reader to  \cite{Fuchs}). 
As a result one  has a Verlinde formula \cite{erik}
for integer fusion coefficients
in terms of (generally non-integer) matrix elements of the 
Verlinde matrix $S$.

A central goal   of the present paper is to understand  the role  
of the group $SL(2,\ZZ)$ in the representation theory of $D(G)$
without reference to conformal field theory. 
Our starting point here is the set of characters of irreducible
representations of $D(G)$. The group $SL(2,\ZZ)$ acts on this 
set in a geometrically natural way.  We identify  two generators 
$S$ and $T$ of this action (satisfying$(ST)^3=S^2$ and $S^4=1$)
which play a natural role in the theory of 
$D(G)$  and its dual  $D(G)^\star$.
It was already noted in 
\cite{DPR} that the  diagonal matrix $T$ is related to 
the central ribbon element of $D(G)$. As a vector spaces, both
$D(G)$ and $D(G)^\star$ can be identified with the space $C(G\times G)$
of $\CC$-valued functions on $G\times G$, and here we show 
that $S$ can be extended  to an  automorphism of the vector space  
$C(G\times G)$.   
We prove a  convolution theorem 
for this extension which shows that it  has a natural interpretation 
as a Fourier transform. 
Returning to the set of characters we show that it
can be given a ring structure
in two dual ways. One, using the algebra multiplication,
is essentially determined by Schur orthogonality relations.
The other, using the dual multiplication, encodes the fusion 
rules of $D(G)$. Our Fourier transform relates the two ring
structures, and we use this  to give a very short proof of the 
Verlinde formula for $D(G)$.

Quantum doubles can also be defined for locally compact
groups $G$ \cite{KM} and we have used a  notation which 
anticipates the generalisation of the arguments given here
to quantum doubles of locally compact groups. There are 
a number of technical problems, however, which we intend
 to address in a future  publication.  
Finally, we observe that  the Fourier transform we will
define in this paper is related to    the  non-abelian Fourier transform 
 defined  by Lusztig  in the context of finite 
group theory, see \cite{Lusztig} and \cite{charred},  and to the  
quantum Fourier transform   defined  by Lyubashenko  in \cite{Ly2}
and discussed further by Lyubashenko and Majid  in
\cite{LyMa2} and \cite{LyMa1}.
We will clarify the relationship between these definitions and 
ours  in the course of the paper.
There are several other places in the  literature  where
Fourier transforms are  discussed in the context of quantum groups.
The focus of the papers \cite{KeMa} and \cite{BKLT} is braided quantum groups
and thus different from ours.
In section 3.4  of  \cite{Chryss} a Fourier transform is defined for finite
dimensional Hopf algebras. However, when applied to the quantum double of a
finite group that definition yields something essentially different
from our  Fourier transform.

\section{The quantum double of a finite group}

Let $G$ be a finite group, with
invariant measure  given by
\bea
\label{maat}
\int_G f(z)\,dz:=|G|^{-1}\,\sum_{z\in G} f(z).
\eea
We will use   delta   functions $\de_x$ on $G$,
normalised so that   $\de_x(y)=0$ if $x\ne y$ and
$\de_x(x)=|G|$.

The quantum double $D(G)$ of a finite group $G$ was first 
studied in detail  in \cite{DPR}.
The definitions we are about to give are  equivalent to the ones 
given there, but 
we adopt a different notation. The advantage of our notation is 
that it easily generalises to the case where $G$ is a locally
compact Lie group \cite{KM}. 
As a linear  space we identify the quantum double $D(G)$ of $G$
with $C(G\times G)$. On $D(G)$ we have a non-degenerate pairing
\bea
\lan f_1,f_2 \ran:=\int_G\!\int_G f_1(x,y)\,f_2(x,y)\,dx\,dy.
\eea
By this pairing we can also identify the dual $D(G)^\star$ of $D(G)$ with
$C(G\times G)$ as a linear space.

On $D(G)$ we have multiplication $\bullet$, identity 1, 
comultiplication $\De$,
counit $\ep$, antipode $\ka$ and involution ${}^*$ given by:
\bea 
\label{algebra}
(f_1\bullet f_2)(x,y)&:=&\int_G f_1(x,z)\,f_2(z^{-1}xz,z^{-1}y)\,dz, 
\nonumber \\
1(x,y)&:=&\de_e(y),\nonumber \\
(\De f)(x_1,y_1;x_2,y_2)&:=&f(x_1x_2,y_1)\,\de_{y_1}(y_2).\nonumber \\
\ep(f)&:=&\int_G f(e,y)\,dy,\nonumber \\
(\ka f)(x,y)&:=&f(y^{-1}x^{-1}y,y^{-1}),\nonumber \\
f^*(x,y)&:=&\wb{f(y^{-1}xy,y^{-1})}.
\eea
By duality we have multiplication $\star$, identity $\iota$, comultiplication
$\De^\star$, counit $\ep^\star$, antipode $\ka^\star$ and involution ${}^\circ$
on $D(G)^\star$:
\bea
\label{coalgebra}
(f_1\star f_2)(x,y)&:=&\int_G f_1(z,y)\,f_2(z^{-1}x,y)\,dz,\nonumber \\
\iota(x,y)&:=&\de_e(x),\nonumber \\
(\De^\star f)(x_1,y_1;x_2,y_2)&:=&f(x_1,y_1y_2)\,\de_{x_2}(y_1^{-1}x_1y_1),
\nonumber \\
\ep^\star(f)&:=&\int_G f(x,e)\,dx,\nonumber \\
(\ka^\star  f)(x,y)&:=&f(y^{-1}x^{-1}y,y^{-1}).\nonumber \\
f^\circ(x,y)&:=&\wb{f(x^{-1},y)}.
\eea

Later, we will also  refer to  the ribbon algebra structure of 
$D(G)$. Following the conventions  for 
ribbon Hopf algebras of sect.  4.2  in \cite{qugroups},
 we define  the universal
$R$-matrix $R\in D(G)\otimes D(G)$:
\bea
R(x_1,y_1;x_2,y_2) = \de_e(y_1)\de_e(x_1y_2^{-1})
\eea
and the central ribbon element $c \in D(G)$:
\bea
\label{randc}
c(x,y) = \bullet\circ(\kappa\otimes{\rm id})(R_{21}) =  \de_e(xy),
\eea
where $R_{21}(x_1,y_1;x_2,y_2):=R(x_2,y_2;x_1,y_1)$. 
The monodromy element $Q\in D(G)\times D(G)$ is 
\bea
\label{monodromy}
Q(x_1,y_1;x_2,y_2):=\bigl(R_{21}\bullet R\bigr)
 (x_1,y_1;x_2,y_2)=\de_{y_1}(x_2)
\de_{y_2}(x_2^{-1}x_1x_2).
\eea
Together with $c$, it  satisfies the ribbon relation
\bea
\label{ribbon}
\Delta c=Q^{-1} \bullet \bigl(  c\otimes c \bigr).
\eea

In the representation theory of $D(G)$ and $D(G)^\star$
 an important role 
is played by the   Haar functionals  $h^\star \colon D(G)^\star \to\CC$
and $h \colon D(G) \to \CC$, respectively. They are given by
\bea
\label{haar}
h^\star (f):=\int_G f(e,y)\,dy \quad \mbox{and} \quad
h (f):=\int_G f(x,e)\,dx
\eea
Here we have chosen the  normalisation 
$h^\star(\iota)=h(1)=|G|$.
Direct computation shows left and right invariance:
\bea
\label{haarprop}
\bigl((h^\star\ten\id)\circ\De^\star \bigr)(f)=h^\star(f)\,\iota=
\bigl((\id\ten h^\star)\circ\De^\star \bigr)(f),
\eea
and similarly for $h$. Furthermore,
centrality, positivity and faithfulness of $h$ and $h^\star$
follow from the formulae
\bea
\label{inpro}
h(f_1\bullet f_2^*)=h^\star(f_1\star f_2^\circ)=
\int_G\int_G f_1(x,y)\,\wb{f_2(x,y)}\,dx\,dy
\eea
and 
\bea 
h^\star(f\star f^\circ)=\int_G\int_G |f(x,y)|^2\,dx\,dy.
\eea
Thus, $C(G\times G)$ has a hermitian inner product
\bea
\label{inprod}
(f_1,f_2)\mapsto h^\star(f_1\star f_2^\circ)=\lan f_1,\wb{f_2}\ran.
\eea
From the existence of faithful
positive linear functionals on $D(G)$ and $D(G)^\star$
(namely $h$ and $h^\star$) we can conclude that $D(G)$ and
$D(G)^\star$ both have a faithful $*$-representation on a finite dimensional
Hilbert space, so they are $C^*$-algebras. Therefore, the theory of
Woronowicz \cite{woronowicz} 
for compact matrix quantum groups holds both for $D(G)$
and $D(G)^\star$. Moreover, this theory simplifies because we are in the
finite dimensional case, see  \cite{woronowicz}, Appendix A.2,  and
\cite{vanDaele}.
These simplifications are already evident in our explicit results that
$\ka^2=\id$ and $(\ka^\star)^2=\id$, and that $h$ and $h^\star$ are central.
Furthermore, in \cite{vanDaele} after the proof of Proposition 2.2,
van Daele  gives a
formula in terms of dual bases for the Haar functional in the finite
dimensional $C^*$-algebra case.
For $D(G)^\star$ this can be written as
\bea
h^\star(f)={\rm const.} \sum_{x,y\in G}(f\star \de_{x,y})(x,y).
\eea
A simple calculation indeed shows that this agrees with (\ref{haar}).
An analogous formula holds for $h$.

\section{Irreducible representations and their characters}

The irreducible representations of the quantum double of a finite
group were classified in \cite{DPR}. We adopt some of the  notation 
used there, but for our description of the  representations 
we follow the approach  used in the discussion of the double
of a locally compact group in  \cite{KM} and \cite{KBM}.
 Thus let $\{C_{A}\}_{A=0,...,p}$ be 
the  set of conjugacy classes in $G$, with $C_0=\{e\}$. In each
$C_A$ pick an element $g_A$ and write $N_A$ for the centraliser of 
$g_A$ in $G$. For later use it is also convenient to  pick, for 
each $x\in G$,
an element $B_x\in G$ such that  
\bea
\label{bdef}
B_xg_A B_x^{-1} =x \quad \mbox{if}  \quad  x \in C_A.
\eea
Write $q_A$ for the number of irreducible representations of $N_A$ 
 and let $\{\pi_\alpha\}$  be a complete set
of  such representations. The label $\al$ is a positive integer 
running from 1 (for  the trivial representation) to $q_A$.
We denote the carrier spaces by  $V_{\alpha}$ and their  dimensions 
by $d_{\al}$. Then
the irreducible representations $\pi^A_{\alpha}$
of $D(G)$  are labelled by pairs
$(A,\alpha)$ of conjugacy classes and centraliser representations.
The carrier space  $V^A_{\al}$ of  $\pi^A_{\alpha}$ is 
\bea
\label{irrep}
V^A_{\al}:=\{\phi:G\rightarrow V_{\al}|\phi(xn)=\pi_\al(n^{-1})\phi(x),
\forall x\in G,\forall n\in N_A\},
\eea
and the action of an element $f\in D(G)$ is
\bea
\label{action}
\left(\pi^A_{\al}(f)\phi\right)(x):=\int_G f(x g_A x^{-1},z)
\, \phi(z^{-1}x)\, dz.
\eea
The set $\{\pi^A_{\al}\}$ is  a complete set of mutually inequivalent
irreducible matrix $*$-representations of $D(G)$. 
We write $d_{A,\al}$ for the   dimension of $V^A_{\al}$ and note that 
 $d_{A,\al} = |C_A|\cdot d_\al$.
Then, after choosing an orthonormal basis of $V^A_{\al}$,
$\pi^A_\al$  can be represented by a matrix $(\pi^A_\al)_{ij}$,
$i,j=1,...,d_{A,\al}$. The matrix elements $(\pi_{\al}^A)_{ij}$ 
are elements of $D(G)^\star $ and we  write 
 $M_{A,\al}$ for the span of the matrix elements $(\pi^A_\al)_{ij}$
($i,j=1,\ldots,d_{A,\al}$). Then it follows from Woronowicz's general
theory that  $D(G)^\star$ is the orthogonal direct sum
of the spaces $M_{A,\al}$. Finally we   define the 
character
\bea
\label{chardef}
\chi^A_\al=\sum_{i=1}^{d_{A,\al}}( \pi^A_\al)_{ii}.
\eea

Characters play a fundamental role in the following discussion. {}From
 \cite{KBM}  we have the following  formula:
\bea
\label{char}
\chi^A_{\al}(f)= \int_G\int_{N_A}f(zg_Az^{-1}, znz^{-1})
\,\chi_{\al}(n)\,dn\,dz.
\eea 
Changing integration variables, this can be rewritten as 
\bea
\label{chara}
\chi^A_{\al}(f)= \int_G\int_G\, f(v,w) \,
{\bf 1}_A(v)\,\,\de_e(vwv^{-1}w^{-1})\,\,\chi_\al
(B^{-1}_v wB_v) \,dv\,dw,
\eea
where 
${\bf 1}_A$ is the characteristic function of the conjugacy class $C_A$,
normalised so  that 
${\bf 1}_A(v)= 1 $  if $v\in C_A$ and   ${\bf 1}_A(v)= 0 $ otherwise.

By definition, characters are elements of $D^\star (G)$. Using the 
pairing $\langle\, ,\rangle$ we can therefore identify them with functions
on $G\times G$. To do this explicitly we 
 insert a delta function for $f$, 
\bea
\label{dprchar}
\chi^A_{\al}(x,y):=\chi^A_{\al}(\de_x\otimes\de_y) =
 \de_e(xyx^{-1}y^{-1})\,{\bf 1}_A(x) \, \chi_{\al}(B_x^{-1}yB_x),
\eea
and  reproduce the character formula given in \cite{DPR}. 
One checks that  the characters
enjoy the orthogonality relation
\bea
\label{ortho}
\langle\chi_{\alpha}^A,\wb{\chi_{\be}^B}\rangle =
|G|\,\,\delta_{\al\be}\,\de_{AB}.
\eea

As elements of  $D( G)^\star$ characters have the 
property that  they are cocentral, i.e. they satisfy
$\De^\star \chi^A_\al = \sigma \circ \De^\star \chi^A_\al$, 
where $\sigma: D(G)^\star
\times D(G)^\star\rightarrow D(G)^\star \times D(G)^\star  $ is 
the flip operation, $\sigma(\lambda,\mu) = (\mu,\lambda)$.
 Using again the identification of $D(G)^\star$
with $C(G\times G)$, the cocentrality of the  characters (\ref{dprchar})
means that their support lies in 
\bea
\label{comm} 
G_{\rm comm}:=\{(x,y)\in G\times G|xy=yx\}.
\eea
and that they are 
invariant  under the  simultaneous
conjugation of both arguments, in symbols 
$\chi^A_\al(gxg^{-1},gyg^{-1}) =\chi^A_\al(x,y)$ for all $g,x,y\in G$. 
These properties are  also evident in the explicit expression (\ref{dprchar}).
We write $C(G_{\rm comm})$ for the  space of functions in $C(G\times G)$
whose support lies in $G_{\rm comm}$, and  $C^{0}(G_{\rm comm})$
for the  space  of functions in $C(G_{\rm comm})$ which 
are invariant under the  simultaneous conjugation of both arguments.
It follows from the remark in  \cite{woronowicz} after Corollary 5.10  that  
the characters in fact span $C^{0}(G_{\rm comm})$. It is instructive 
to see this explicitly. The dimension of $C^{0}(G_{\rm comm})$
is equal to the number of $G$-conjugacy
classes in $G_{\rm comm}$. To count these,
introduce an integer label $a$  for  $N_A$-conjugacy classes in $N_A$.
Since the number of such conjugacy classes is equal 
the number of irreducible representations, $a$ runs from $1$ to $q_A$.
In the $N_A$-conjugacy class labelled by $a$ pick  an element $g^a_A$.
Then every $G$-conjugacy class in $ G_{\rm comm}$ contains a unique element
of the form $(g_A,g_A^a)$ for suitable $A$ and $a$.
The number of conjugacy classes is  therefore 
\bea
\mbox{dim}\bigr( C^{0}(G_{\rm comm})\bigr) = \sum_{A=0}^p q_A.
\eea
This  is also  the number of irreducible representations $(A,\al)$ of 
$D(G)$ and hence  the number of characters.
 By the orthogonality relation, the characters
 are certainly linearly independent  and therefore form
an orthogonal basis of   $C^{0}(G_{\rm comm})$.

The vector space $C^{0}(G_{\rm comm})$ 
is closed under both the multiplication $\bullet $
and the dual multiplication 
 $\star $. This means that the vector space spanned by
the  characters  can be given two algebra structures which 
are  dual to each other. Both these algebras are initially
defined over the field $\CC$, but in the following section we will show
that  the structure constants for both algebras are integers.
Therefore the integer linear combinations 
of characters form a ring over $\ZZ$ 
under both the multiplications $\bullet$
and $\star$.

\section{Character rings}

Again following the general theory given in \cite{woronowicz},  
 the matrix elements $(\pi_{\al}^A)_{ij}$  form  a complete set
of mutually inequivalent irreducible matrix 
corepresentations $\pi^A_{\al}$ of
$D(G)^\star $. We therefore have
\bea
\label{costar}
\De^\star \bigl((\pi^A_\al)_{ij}\bigr)=\sum_k 
(\pi^A_\al)_{ik}\ten(\pi^A_\al)_{kj}.
\eea
The quantum analogues of Schur's orthogonality relations, given  
in  \cite{woronowicz},   simplify for $D(G)$
because \hbox{$(\ka^\star )^2=\id$}. We have
\bea
\lan (\pi^A_\al)_{ij}, (\wb{\pi}^B_\be)_{kl}\ran=  h^\star
((\pi^A_\al)_{ij} \star 
(\pi^B_\be)_{kl}^\circ)=\de_{\al\be}\,\de_{AB}\,
\de_{ik}\,\de_{jl}\,
h^\star(\iota)/d_{A,\al}.
\eea
Note that  due to a standard theorem in
the theory of finite groups, see  \cite{Serre} or \cite{groupth}, 
the ratio 
\bea
\label{ndef}
n^A_\al:=h^\star(\iota)/d_{A,\al}=|N_A|/d_\al
\eea
is an integer. These  relations 
are sufficient to establish the following theorem.

\begin{theorem}
\label{projthm}
The map $f\mapsto (n^A_\al)^{-1}\,
\chi^A_\al\bullet f$ is the orthogonal projection of $D(G)^\star$ onto 
$M_{A,\al}$.
\end{theorem}
{\bf Proof}\quad We have:
$$
\lan\chi^A_\al\bullet(\pi^B_\be)_{ij},\wb{(\pi^C_\ga)_{kl}}\ran
=
\lan \chi^A_\al\ten(\pi^B_\be)_{ij},\De^\star  \wb{(\pi^C_\ga)_{kl}}\ran
=
\sum_r \lan\chi^A_\al,\wb{(\pi^C_\ga)_{kr}}\ran\,
\lan(\pi^B_\be)_{ij},\wb{(\pi^C_\ga)_{rl}}\ran
$$
$$
=
\sum_{m,r}\lan(\pi^A_\al)_{mm},\wb{(\pi^C_\ga)_{kr}}\ran\,
\lan(\pi^B_\be)_{ij},\wb{(\pi^C_\ga)_{rl}}\ran
=\de_{AC}\,\de_{BC}\,
\de_{\al\ga}\,\de_{\be\ga} \,\de_{ik}\,\de_{jl}\, \bigl(n^A_\al\bigr)^2.
\eqno \halmos
$$

As an immediate consequence  we note:  
\begin{lemma}
\label{bul}
The  characters of the quantum double $D(G)$ of a finite group $G$ 
form a ring over $\ZZ$ under the multiplication
$\bullet$. The multiplication rule is
\bea
\label{bulalg}
\chi^A_\al\bullet \chi^B_\be = \de_{AB}\,\de_{\al\be}\,
n^A_\al \,\chi^A_\al.
\eea
\end{lemma}

Next consider  the  algebra structure of the characters under the 
dual multiplication $\star $. This is related to the tensor product 
decomposition into irreducible representations:
\bea
\label{fusion}
\pi^A_\al \otimes \pi^B_\be \simeq \bigoplus_{C,\gamma}
N_{\al\be C}^{A B \ga} \pi^C_\ga.
\eea
We will refer to
 the non-negative integers  $N_{\al\be C}^{A B \ga}$ as  fusion
coefficients.
By definition of the dual multiplication we have, 
for  $\pi,\rho \in D(G)^\star $
and $f\in D(G)$
\bea
\label{reminder}
\lan\pi\otimes\rho, \De (f)\ran = \lan\pi\star \rho, f\ran.
\eea
Thus, upon taking the trace we find that for all $f\in D(G)$
\bea
\tr\bigl(\pi^A_\al \otimes \pi^B_\be \bigl(\De(f)\bigr)\bigr)=
\sum_{C,\gamma} \,N_{\al\be C}^{A B \ga}\, \tr\bigl( \pi^C_\ga (f)\bigr).
\eea
Using (\ref{reminder})  we deduce the following  lemma. 
\begin{lemma} 
\label{starr}
The  characters of the quantum double $D(G)$ of a finite group $G$ 
form a ring over $\ZZ$ under the multiplication
$\star$. The multiplication rule is
 determined by the fusion coefficients
$N_{\al\be C}^{A B \ga}$:
\bea
\label{starmult}
\chi^A_\al \star \chi^B_\be = \sum_{C,\gamma} \,N_{\al\be C}^{A B \ga}
\,\chi^C_\ga.
\eea
\end{lemma}

\section{$SL(2,\ZZ)$-action, Fourier transform and the Verlinde formula}

In this section, a central role is played 
by a  natural action of the group  $SL(2,\ZZ)$  
of  integer unimodular $2\times 2$ matrices on  space  $C(G_{\rm comm})$.
Let 
\bea 
 M= \left(
\begin{array}{cc}
a&b\\
c&d 
\end{array}\right),
\eea
with $a,b,c,d$ integers such that $ad-bc=1$,
be a generic element of $SL(2,\ZZ)$  and define a  right action
on $G_{\rm comm}$ via
\bea
(x,y) M: = (x^ay^c,x^b y^d).
\eea 
This induces an action of $M \in SL(2,\ZZ)$ on elements 
 $f\in C(G_{\rm comm})$, which we write as 
\bea 
(Mf)(x,y):= f(x^ay^c,x^b y^d).
\eea
The generators 
\bea
 S=
\left(
\begin{array}{cc}
0&-1\\
1&0
\end{array}\right), ~~~~~
 T=
\left(
\begin{array}{cc}
1&1\\
0&1
\end{array}\right)
\eea
satisfy the modular relation
\bea
\label{modrel}
( ST)^3= S^2
\eea
 and $S^4=1$, 
and both have a natural interpretation
within the quantum double $D(G)$. To see this, first note that
the actions 
\bea 
\label{four}
(Sf)(x,y) =f(y,x^{-1})
\eea
and 
\bea
\label{spinf}
(Tf)(x,y)= f(x,xy)
\eea
also make sense for general $f\in C(G\times G)$.
We keep  the notation (\ref{four}) and (\ref{spinf}) 
even when the arguments $x$ and $y$ do not commute. Note that
$S$ and $T$ are unitary operators on $C(G\times G)$ with the inner
product (\ref{inprod}).
Moreover, one  finds that the action of $T$ on $f\in C(G\times G)$ is equal to 
the multiplication of $f$ by the central element $c$ (\ref{randc}):
\bea  
\label{spinff}
Tf =c\bullet f.
\eea
Acting on $ C(G\times G)$, $S$ and $T$ no longer satisfy the 
modular relation (\ref{modrel})  but we still have $S^4=1$. 
This last property and  the following convolution
theorem  are our reasons for calling $S$ a Fourier transform.

\begin{theorem}
\label{invol}
If $f,g\in C(G\times G)$ and ${\rm supp}(f)\in G_{\rm comm}$ then
\bea 
\label{fourthe}  
S(f\star g)=S(f) \bullet S( g)\quad \hbox{and } \quad
S(f\bullet g)=S(g)\star S( f).
\eea
\end{theorem}

\mLP
{\bf Proof}\quad  If the first factor in a $\bullet$-product has
support in  $ G_{\rm comm}$, the formulae simplify, yielding 
\bea
 \bigl(S(f)\bullet   S(g)\bigr)(x,y)
&=&\int_G f(z,x^{-1}) g(z^{-1}y, x^{-1}) \,dz\nonumber \\
&=& \bigl( S(f\star g)\bigr)(x,y). \nonumber
\eea
Similarly we have 
\bea
\bigl(S(g)\star S(f) \bigr)
(x,y) &=& \int_G g(y,z^{-1}) f(y,x^{-1}z)\,dz \nonumber \\
&=&  \int_G f( y, w ) \,g (y,w^{-1} x^{-1})\,dw \nonumber \\
&=& \bigl(S(f\bullet  g)\bigr)(x,y),\nonumber
\eea
where we have again used 
${\rm supp}(f)\in  G_{\rm comm}$ and exploited the 
invariance of the measure under
the change of integration  variable  $w= x^{-1}z$.\hhalmos

Since $S$ leaves the space of  functions with support in $G_{\rm comm}$
invariant, we deduce the following 
\begin{corollary}
\label{sinv}
If $f,g\in C(G\times G)$ and ${\rm supp}(f)\in G_{\rm comm}$ then
\bea
S^{-1}(f\bullet g)=S^{-1}(f)\star S^{-1}(g)\quad \hbox{and } \quad 
S^{-1}(g\star f)= S^{-1}(f)\bullet S^{-1}(g).
\eea
\end{corollary}

At this point it is instructive to make contact with 
related discussions of Fourier transforms  in the literature.
By defining a slight variant of the operator $S$ we can establish 
a connection  with
the non-abelian Fourier transform given  by
Lusztig  in   \cite{Lusztig} and \cite{charred}.
For $f\in C(G\times G)$ put
\bea
(Uf)(x,y):=f(y,x),\quad
(J_1f)(x,y):=f(x^{-1},y),\quad
(J_2f)(x,y):=f(x,y^{-1}).
\eea
Then $U=J_2S=SJ_1$. The operators $U$, $J_1$, $J_2$ correspond to
$2\times 2$ matrices with integer entries but with determinant $-1$.
 Note also that the operators
$U$, $J_1$ and $J_2$, like $S$, leave the
space of functions
with support in $G_{\rm comm}$ invariant and commute with conjugations by
elements of $G$. We can therefore in particular 
consider the restriction of $U$
to $C^0(G_{\rm comm})$. 
It turns out that the matrix elements of $U$  with respect to  the basis of 
characters formally coincide with
Lusztig's Fourier kernel:
\bea
\label{umatrix}
U^{BA}_{\be\al}:&=&|G|^{-1}\lan U\chi_\al^A,\wb{\chi_\be^B}\ran
\nonumber \\
&=&{1\over |N_A|\,|N_B|}\,
\sum_{\textstyle{g\in G\atop g_A\,gg_Bg^{-1}=gg_Bg^{-1}\,g_A}}
\chi_\al(gg_Bg^{-1})\,\wb{\chi_\be(g^{-1}g_Ag)} \nonumber \\
&=&: \{(g_A,\al),(g_B,\be)\}.
\eea
However, Lusztig takes $\{\,,\,\}$ with values in
the field $\wb \QQ_l$, i.e., in an algebraic completion of the field
$\QQ_l$ of l-adic numbers. He  also has  a bar operation on $\wb \QQ_l$.
 
Straightforward computations shows that 
the following analogues of Theorem (\ref{invol}) hold.  For 
 $f,g\in C(G\times G)$ we have 
\bea
J_1(f)\star J_1(g)=J_1(g\star f),\quad
J_1(f)\bullet J_1(g)=J_1(f\bullet g)\quad
\mbox{and}\quad J_2(f)\star J_2(g)=J_2(f\star g).
\eea
If supp$(f)\in G_{\rm comm}$ one checks that furthermore
\bea
U(f)\star U(g)=U(f\bullet g)\quad{\rm and}\quad
U(f)\bullet U(g)=U(f\star g).
\eea
Finally,  if supp$(f)$ and supp$(g)\in G_{\rm comm}$ then  
\bea
J_2(f)\bullet J_2(g)=J_2(g\bullet f).
\eea

In the abstract setting of tensor categories, a quantum Fourier
transform 
 was defined by Lyubashenko in \cite{Ly2} and discussed further
 in \cite{LyMa2} and \cite{LyMa1}
for finite-dimensional 
factorisable ribbon Hopf algebras.  Applied to $D(G)$ and  in our
notation their formula for the Fourier transform $\tilde S$
of an element $f$ of $D(G)$ reads
\bea
\label{shahn}
\tilde Sf: = (1\otimes h )\circ \bigl(R^{-1}\bullet (1\otimes f) \bullet
R_{21}^{-1}\bigr).
\eea
An explicit calculation  shows that 
\bea 
\label{fourti}
\bigl(\tilde S f\bigr)(x,y) = f(xy^{-1}x^{-1}, x).
\eea
and therefore the relation between $\tilde S$ and $S$ can
be expressed via
\bea
\label{shrel}
\tilde S = \kappa \circ S.
\eea
The fourth power of the  map $\tilde S$  is not equal to the identity,
but if, following \cite{LyMa2}, one defines $\tilde T = T^{-1}$ one
finds that the modular relation $ (\tilde S \tilde  T)^3 = {\tilde S}^2$
is satisfied. 
Convolution formulae  similar to ours 
can be proven for the map (\ref{shahn}). Again, at least one of the two
elements $f$ and $g$ to be multiplied has to have support 
in $ G_{\rm comm}$. Finally, we observe that restricted to $C(G_{\rm comm})$,
the map $\tilde S$ agrees with our $S^{-1}$.

For the rest of this paper we focus our attention on
the space $C^{0}(G_{\rm comm})$. 
In particular we consider the restriction of the map 
$S$ to $C^{0}(G_{\rm comm})$, and again denote it by $S$.
The characters
form a natural orthogonal basis of the  space 
$C^{0}(G_{\rm comm})$,
 and we
define the matrix $S^{BA}_{\be\al}$ as the 
matrix representing the map $S$ on the  basis of characters:
\bea
\label{verlmatrix}
S^{BA}_{\be\al} := |G|^{-1}\,\lan S\chi^A_\al, \wb{\chi^B_\be}\ran.
\eea 
Here the normalisation is chosen so that 
$S\chi^A_\al=\sum_{B,\beta} S^{B A}_{\be\al} \chi^B_\beta$.
The matrix $S^{BA}_{\be\al}$ is unitary because the map $S$ is.
Using the explicit expression for the characters (\ref{dprchar})
one finds the following  formula, first  given in \cite{DPR}:
\bea
\label{matform}
S^{BA}_{\be\al}=\int_G\int _G \de_e(xyx^{-1}y^{-1})\,
\, {\bf 1 }_A(x)\, {\bf 1}_B(y)\,
\,\wb{\chi}_{\al}(B_x^{-1}yB_x)\,\,\wb{ \chi}_\be(B_y^{-1}x B_y)\,\,dx\,dy.
\eea
This expression shows that the matrix $S^{AB}_{\al\be}$  is symmetric,
$S^{AB}_{\al\be}= S^{B A}_{\be\al }$. Since it is also unitary
its inverse is given by its complex conjugate. We can also read off 
the useful relation
\bea 
\label{dimform}  
S^{ 0A}_{1 \al}=\frac{1}{n^A_\al}.
\eea

 Armed with this notation, we can now use the convolution
theorem to relate the $\bullet$ and $\star $-ring structures 
of the characters. The result is the  Verlinde formula.

\begin{theorem}[Verlinde Formula] \quad 
Acting on characters, the inverse  Fourier transform 
 $S^{-1}$
diagonalises the fusion rules of $D(G)$. 
The fusion coefficients can  be expressed in terms
of the  matrix $S^{AB}_{\al\be}$ :
\bea
N^{AB\ga}_{\al\be C}=
\sum_{D,\de}\frac{ S^{DA}_{\de\al}\,S^{DB}_{\de\be}\,
\wb{S}^{CD}_{\ga\de}}{S^{0D}_{1\de}}.
\eea
\end{theorem}

\mLP
{\bf Proof} \quad
It follows from  the definition of $S$ and from Lemma (\ref{bul}) that 
$$
 S\chi^A_\al\bullet \chi^B_\be =\frac{ S^{BA}_{\be\al}} {S^{0B}_{1\be}}\,\,
\chi^B_\be.
$$ 
Now apply  $S^{-1}$ to both sides   and use 
 the first  formula in Corollary (\ref{sinv})
to obtain 
$$
\chi^A_\al\star S^{-1}\chi^B_\be =\frac{ S^{AB}_{\al\be}} {S^{0B}_{1\be}}\,\,
S^{-1}\chi^B_\be,
$$
yielding the diagonalised fusion rules, with eigenvalues
$ S^{AB}_{\al\be}/ S^{0B}_{1\be}$. A quick derivation
of the formula for 
the fusion coefficients  follows again 
 from  the definition of $S$ and from Lemma (\ref{bul}):
$$
S\chi^A_\al\bullet S\chi^B_\be   = 
\sum_{D,\de}\frac{ S^{DA}_{\de\al}\, S^{DB}_{\de\be}} {S^{0D}_{1\de}}\,\,
\chi^D_\de.
$$
Again apply   $S^{-1}$ to both sides   and use 
 the first  formula in Corollary (\ref{sinv})
to obtain 
$$
\chi^A_\al\star\chi^B_\be = \sum_{C\ga} \,\left(
\sum_{D,\de}\frac{ S^{DA}_{\de\al}\,S^{DB}_{\de\be}\,
\wb{S}^{CD}_{\ga\de}}{S^{0D}_{1\de}}
\right) \,\chi_\ga^C.
$$
Comparing  this expression 
with  (\ref{starmult}) shows that the expression
in brackets is equal to the fusion coefficient 
$N^{AB\ga}_{\al\be C}$.\hhalmos

\noindent{\bf Remark} \quad 
There is an interesting connection with Lusztig's matrix 
$U^{AB}_{\al \be}$ (\ref{umatrix})
here. We find that $\wb{U_{\al\be}^{AB}}=U_{\be\al}^{BA}= 
(U^{-1})_{\be\al}^{BA}$ and that
$U_{\al\be}^{AB}=S_{\al \wb \be}^{AB}$.
Verlinde's formula can also be expressed in terms of $U$:
\bea
\chi_\al^A\star U\chi_\be^B={U_{\be\al}^{BA}\over U_{ 1\be}^{0B}}\,
U\chi_\be^B
\eea
and
\bea
N_{\al\be C}^{AB\ga}=
\sum_{D,\de}{{U_{\de\al}^{DA}}\,\,{U_{\de\be}^{DB}}\,\,
U_{\ga\de}^{CD}\over
U_{1\de}^{0D}}\,.
\eea

\vspace{1.3cm}

The simplicity of our proof of the Verlinde formula 
shows that the Fourier transform
$S$ is a very natural   tool for proving the Verlinde formula  for 
$D(G)$. While we have restricted attention to a particular 
ribbon Hopf algebra here, we  have tried to indicate as far
as possible how  our definitions and equations for $D(G)$
can be formulated using only natural operations (such the 
Haar measure, the antipode, the universal $R$-matrix
or the central ribbon element) 
which exist  for a large class of (quasi-triangular ribbon)
Hopf algebras. More generally it is natural to ask
for which class of 
(quasi) Hopf algebras a Fourier transform with analogous 
properties  can be defined. In view of the tight connection
between fusion rules in rational conformal field theory and 
tensor decomposition rules in (quasi) Hopf algebras (see e.g. 
\cite{FGV}, or  \cite{Fuchs} for a review)  such a  generalised Fourier
transform, if it exists, could be expected to  play an important role in
both Hopf algebra theory and conformal field theory.

\vspace{1cm}

\noindent {\bf Acknowledgements}
 
\noindent BJS was a post doc at the Instituut voor Theoretische Fysica
of the University of Amsterdam 
while the research for this paper was carried out. 
BJS thanks J\"urgen Fuchs for illuminating discussions about 
fusion rules.



\end{document}